\documentclass[11pt]{article}

\usepackage{graphicx}

\newcommand{\bea}{\begin{eqnarray}}
\newcommand{\eea}{\end{eqnarray}}
\newcommand{\be}{\begin{eqpation}}
\newcommand{\ee}{\end{equation}}
\newcommand{\tras}{\mbox{\sc t}}

\newcommand\bu {{\bf{u}}}

\newcommand\bx {{\bf{x}}}

\newcommand\bX {{\bf{x}}}

\newcommand\wf {\widehat{f}}

\newcommand\wg {\widehat{g}}

\newcommand\bphi {\mbox{\boldmath $\phi$}}

\newcommand\bbeta {\mbox{\boldmath $\beta$}}

\newcommand\weta {\widehat{\eta}}

\newcommand\wbeta {\widehat{\mbox{\boldmath $\beta$}}}

\newcommand\wbphi {\widehat{\bphi}}

\newcommand{\wfi} {\widehat{\phi}}

\newcommand{\wbfi} {\widehat{\mbox{\boldmath $\phi$}}}
\newcommand\wtX {\widetilde{\bX}}
\newcommand\wgama {\widehat{\gamma}}
\newcommand\wtgama {\widetilde{\gamma}}

\newcommand\wtfi {\widetilde{\phi}}

\def\var{\mathop{\rm var}}
\def\MSE{\mathop{\rm MSE}}
\def\real{\hbox{$\displaystyle I\hskip -3pt R$}}

\def\realito{\tiny{\real}}

\newcommand{\convprob}{ \buildrel{p}\over\longrightarrow}
\newcommand{\convdist}{ \buildrel{{\cal D}}\over\longrightarrow}

\def\dst{\displaystyle}
\def\noi{\noindent}

\def\square{\ifmmode\sqr\else{$\sqr$}\fi}
\def\sqr{\vcenter{
         \hrule height.1mm
         \hbox{\vrule width.1mm height2.2mm\kern2.18mm
\vrule width.1mm}
         \hrule height.1mm}}

\usepackage{color}
\begin{document}
\title{Partially linear models on Riemannian manifolds}

\author{Wenceslao Gonzalez-Manteiga$^1$, Guillermo Henry$^2$ and Daniela Rodriguez$^2$ \\
{\small $^1$\sl Universidad de Santiago de Compostela, Spain}\\
{\small $^2$ \sl Facultad de Ciencias Exactas y Naturales, Universidad de Buenos Aires and
CONICET, Argentina}}

\date{}
\maketitle

\begin{abstract}
 In partially linear models the dependence of the  response $y$  on  $(\bx^{\tras},t)$  is modeled through the relationship
 $ y=\bx^{\tras} \bbeta+g(t)+\varepsilon$  where $\varepsilon$ is independent of $(\bx^{\tras},t)$.  In
this paper, estimators of $\bbeta$ and $g$   are  constructed when the explanatory variables $t$ take values on a Riemannian manifold.
Our proposal combine the flexibility of these models with the complex structure of a set of explanatory variables.
We prove that the resulting estimator of $\bbeta$ is  asymptotically normal under the suitable conditions. Through a simulation study, we explored the performance of the estimators. Finally, we applied the studied model to an example based on real dataset.
\end{abstract}

\noindent{\em Key words and phrases:}
Nonparametric estimation, Partly linear models,  Riemannian manifolds.

\section{Introduction}
The partially linear models was introduced by \cite{eng} to analyzed the relationship between the electricity usage and average daily temperature. In recent years, this model has gained a lot of attention in order  to explore the nature of complex nonlinear phenomena. This model  has been widely studied in the literature see for example \cite{sp}, \cite{chen}, \cite{aq} among others. The partially linear models allow modeling the  response variable with a set of predictors that  enter  linearly in the model  while one of them is considered in the model nonparametrically.

 However, in many applications, the predictors variables take values on a Riemannian manifold more than on Euclidean space and  this structure of the variables needs to be taken into account in the estimation procedure. Some examples could be found in meteorology, astronomy, geology and other fields, that include distributions on spheres,  tangent bundles, Lie groups, etc. Research on the statistical analysis of variables  with some one of this structures  was studied by \cite{BP}, \cite{Mardia} and more recently by \cite{HL}, \cite{pennec},\cite{bp} and \cite{hrnp}.

The aim of this work is to study the partially linear models when the explanatory variable $t$ takes values on a Riemannian manifold, i.e.  when the variable to be modeled in a nonparametric way is in a manifold. Our proposal combine the flexibility for these models with the complex structure of a set of explanatory variables.

This paper is organized as follows. In Section \ref{estimadores}, we construct estimates for this models and give a brief summary of the nonparametric estimation on Riemannian manifolds proposed in \cite{bp}. In Section \ref{asymp}, we  present the asymptotic distribution of the regression parameter under regular assumptions on the bandwidth sequence. In Section \ref{simulacion}, we explored the performance of the estimators with a simulation study and we show an example using real data. Also, we review a cross validation procedure for partial linear models. Proofs are given in the Appendix.

\section{Estimators} \label{estimadores}
\subsection{Model and estimators}
Let $(y_i,\bx_i^{\tras},t_i)$ be an i.i.d. random vectors valued in $\real^{p+1} \times M$ with identically  distribution to $(y,\bx^{\tras},t)$, where $(M,g)$ is a Riemannian manifolds of dimension $d$. The partially linear model assume that the relation between the response variable $y_i$ and the covariates $(\bx_i^{\tras},t_i)$ can be represented  as
\begin{equation}
y_i=\bx^{\tras}_i \bbeta+g(t_i)+\varepsilon_i \quad\quad 1\leq i \leq n\;,
\label{semipara}
\end{equation}
where the errors $\varepsilon_{i}$ are independent and independent of $(\bx_i^{\tras},t_i)^{\tras}$, also $E(\varepsilon_i|\bx_i,t_i)=0$.  In many situations, it seems reasonable to suppose that a relationship between the covariates $\bx$ and $t$ exists, so as in \cite{sp}
 and \cite{aq}, we will assume that for $1\leq j \leq p$
\begin{equation}
x_{ij}= \phi_j(t_i)+\eta_{ij}  \quad \quad 1 \le i\le n \, \label{semi2}
\end{equation}
where the errors $\eta_{ij}$ are independent. Denote $\phi_0(\tau)=E(y|t=\tau)$ and $\bphi(t)=(\phi_1(t),\dots,\phi_p(t))$,  then we have that $g(t)=\phi_0(t)-\bphi(t)^{\tras}\bbeta$ and hence, $y-\phi_0(t)=(\bx-\bphi(t))^{\tras}\bbeta+\varepsilon$. This equation suggest estimate  the  unknown functions and parameters  as follows.  Let $\wfi_j(t)$ be the nonparametric estimators of $\phi_j$ for $0\leq j\leq p$. Note that the regression functions correspond to  predictors taking values in a Riemannian manifold,  nonparametric kernel type  estimators adapted to this structure was considered in \cite{bp} and also studied in \cite{hriv}. An overview of this estimators can be found in the following Subsection.

Returned to the estimation of the parameter $\bbeta$, note that  using the  nonparametric estimators of the functions $\phi_j$,  the regression parameter can be estimate considering the least square estimators obtained minimizing
\begin{eqnarray*}
\widehat{\bbeta}=\mbox{arg} \min_{\bbeta} \sum_{i=1}^n[(y_i-\wfi_0(t_i))-(\bX_i-\wbfi(t_i))^{\tras}\bbeta]^2.
\end{eqnarray*}
where $\wbphi(t)=(\wfi_1(t),\dots,\wfi_p(t))$. Then the function $g$ can be estimated as $\wg(t)=\wfi_0(t)-\wbphi(t)^{\tras}\widehat{\bbeta}$. This procedure is consistent with the respective estimators when the explicative variable $t$ take values on Euclidean spaces, i.e. the proposed estimators  reduce to  know estimators introduced by \cite{eng}.

\subsection{Review of Nonparametric estimators on Riemannian manifolds}
\subsubsection{Preliminaries}\label{preli}
As in \cite{hrnp} we consider $({M},g)$ a $d-$dimensional oriented Riemannian manifold without boundary, complete and with positive injectivity radius ($inj_g M>0$ ). From now on, $d_g$ will denote the distance function induced by the metric $g$. Throughout this note, we will consider  the concept of volume density function. For a rigorous definition of this function see \cite{Besse} or \cite{hriv}. If we consider the exponential normal chart $(U,\psi)$ of $(M,g)$ induced by an orthonormal basis
$\{v_1,\dots,v_d\}$ of $T_sM$,  then
$\theta_s(t)={\left|\det g_t\left({\partial}/{\partial
\psi_i}\Big|_t,{\partial}/{\partial \psi_j}\Big|_t
\right)\right|}^{\frac 12}\ ,$
where ${\partial}/{\partial\psi_i}|_t=D_{\alpha_i(0)}exp_s(\dot{\alpha_i}(0))$ with $\alpha_i(u)=exp_s^{-1}(t)+uv_i$ for $t\in U$. For example,  when $M$ is $\real^d$ with the canonical metric, then $\theta_s(t)=1$ for all $s,t\in\real^d$ and also in the case of the cylinder $\theta_s(t)=1$. In \cite{hrnp},  we calculate the volume density on the sphere, in this case, $\theta_s(t)={|sen(d_g(s,t))|}/{d_g(s,t)}$  for $ t\neq s, -s.$
and $\theta_s(\pm s)=1$. See also, \cite{hrnp} for a discussion on the geometric definitions.

\subsubsection{The nonparametric estimators}

Let $(y_1,t_1),\cdots, (y_n,t_n)$ be i.i.d  random objects that take values on $\real\times M$. In order to estimate $r(\tau)=E(y|t=\tau)$, Pelletier \cite{bp} proposed a nonparametric kernel type  estimators. The Pelletier´s idea was to build an analogue of a kernel on $(M,g)$, by using a positive function of $d_g$ distance normalized by the volume density function of $(M,g)$, to take into account the curvature of the manifolds.
More precisely, the nonparametric estimator can be defined as,
\begin{equation}\label{estimadornop}
r_n(t)=\sum_{i=1}^n w_{n,h}(t,t_i) y_i
\end{equation}
with $w_{n,h}(t,t_i)={{\theta^{-1}_t(t_i)}K(d_g(t,t_i)/h)}/[{\sum_{k=1}^n {\theta^{-1}_t(t_k)}K(d_g(t,t_k)/h)}]^{-1}$
where $K:\real \to \real$ is a non-negative function, $\theta_t(s)$  the volume density
function on $(M,g)$ and the bandwidth $h$  is a sequence of real positive numbers such that $\lim_{n\to \infty}h=0$ and $h< inj_g M$, for all $n$.  This last requirement on the bandwidth guarantees that (\ref{estimadornop}) is defined for all $t\in M$.
In \cite{bp}, is  derived an expression for the asymptotic pointwise bias and variance as well as an expression for the asymptotic integrated mean square error. On the other hand, in \cite{hrnp} is proposed a robust version  that generalized
these estimators  and it is obtained  the uniform almost sure consistency over compact set and derived the asymptotic distribution.

\section{Asymptotic behavior}\label{asymp}
The theorem of this section studies the asymptotic behavior of the regression parameter estimator of the model under the following  conditions.

\begin{enumerate}
\item[$H1.$] Let ${M}_0$ be a compact set on ${M}$ such that: $f$ is a bounded function such that $\inf_{t\in M_0}f(t)=A>0$ and
 $\inf_{t,s\in M_0}\theta_t(s)=B>0.$
\item[$H2.$] The sequence $h$ is such that $nh^4\to 0$ and ${nh_n^d}/{\log n}\to \infty $ as $n\to \infty$.
\item[$H3.$] $K: \real \to \real$ is a bounded
nonnegative Lipschitz function of order one, with compact support $[0,1]$   satisfying: $\int_{\realito^d}
K(\|\bu\|)d\bu=1$, $\int_{\realito^d} \bu
K(\|\bu\|)d\bu=\bf{0}$ and $0<\int_{\realito^d}
\|\bu\|^2K(\|\bu\|)d\bu<\infty$.
\item[$H4.$] For any open set $U_0$ of $M_0$ such that $M_0\subset U_0$, the functions $g, \phi_j$ for $1\leq j\leq p$  are of  class $C^2$ on $U_0$.
\item[$H5.$] The errors $\varepsilon_i$ and $\eta_{ij}$   for $1\leq i\leq n$ and $1\leq j \leq p$ are independent and $E|\varepsilon_1|^r+\sum_{j=1}^pE|\eta_{1j}|^r<\infty$ for $r\geq 3$, $\sigma_{\varepsilon}^2=\var(\varepsilon_1)>0$ and $\Sigma=E(\eta_1^{\tras}\eta_1)$ is a positive defined matrix.
\end{enumerate}

\noi \textbf{Remark \ref{asymp}.1.} The fact that $\theta_t(t)=1$ for all $p\in M$ guarantees that the bonded of $\theta$ in $H1$  holds. The assumptions $H2$ and $H3$  are standard assumptions when dealing kernel estimators.

\vspace{0.1cm}
\noindent
\noi \textbf{Theorem \ref{asymp}.1.} Under $H1$ to $H5$ we have that $\sqrt{n}(\wbeta-\bbeta)\convdist N(0,\sigma^2_{\varepsilon}\Sigma^{-1})$. \rm
%\vspace{0.2cm}

\noi \textbf{Remark \ref{asymp}.2.} Note that this theorem is consistent with the respective results in the Euclidean case. The obtained asymptotic distribution can be used to construct a Wald-type statistics to make inference on the regression parameter, that is, when we want to test $H_0: \bbeta=\bbeta_0.$

\section{Real example and Monte Carlo study}\label{simulacion}

\subsection{Selection of the smoothing parameter}\label{cv}

An important issue in any smoothing procedure is the choice of the smoothing parameter. Under a nonparametric regression model with carriers in an Euclidean space, i.e., when $M$ is $\real^d$ with the canonical metric, two commonly used approaches are $L^2$ cross--validation and plug--in methods.  In this section,  we included a cross-validation method  for the choice of the bandwidth in the case of partially linear models. The asymptotic properties of data--driven estimators require further careful investigation and are beyond the scope of this paper.

The cross-validation method constructs an asymptotically optimal data-driven bandwidth, and thus adaptive data-driven estimators, by minimizing
$CV(h) =\sum_{i=1}^n[(y_i-\wfi_{0,-i,h}(t_i))-(\bX_i-\wbfi_{-i,h}(t_i))^{\tras}\widetilde{\bbeta}]^2,$
where $\wfi_{0,-i,h}(t)$ and $\wbfi_{-i,h}(t)=(\wfi_{1,-i,h}(t),\dots,\wfi_{p,-i,h}(t))$ denote the nonparametric estimators computed with bandwidth $h$ using all the data expect the $i-$th observation and $\widetilde{\bbeta}$
minimize $\sum_{i=1}^n[(y_i-\wfi_{0,-i,h}(t_i))-(\bX_i-\wbfi_{-i,h}(t_i))^{\tras}{\bbeta}]^2$ in $\bbeta$.

\subsection{Simulation study}\label{simu}

To evaluate the performance of the estimation procedure, we conduct a simulation study. We consider two models in two different Riemannian manifolds, the sphere and the cylinder endowed with the metric induced by the canonical metric of $\real^3$. We performed 1000 replications of independent samples of size $n=200$ according to the following models:

\noi\bf Sphere case: \rm   The variables $(y_i,x_i,t_i)$ for $1\leq i\leq n$ were generated as
$$
y_i=\beta\;x_i+\exp{\{-(t_{i1}+2t_{i2}+t_{i3})^2\}}+\varepsilon_i \quad \mbox{ and } \quad x_i=t_{i1}+t_{i2}+t_{i3}+\eta_i
$$
where  $t_i=(\cos(\theta_i)\cos(\gamma_i),\sin(\theta_i)\cos(\gamma_i),\sin(\gamma_i))$ with $\theta_i$ and $\gamma_i$ follow a von Mises distribution with means $0$ and $\pi$ and concentration parameters $3$ and $5$, respectively.

\noi \bf Cylinder case: \rm The variables $(y_i,x_i,t_i)$ for $1\leq i\leq n$ were generated as
$$
y_i=\beta\;x_i+s_i^2+\sin(\theta_i)+\varepsilon_i \quad \mbox{ and } \quad x_i=\exp(\theta_i)+\eta_i
$$
where $t_i=(\cos(\theta_i),\sin(\theta_i),s_i)$ with the variables $\theta_i$ follow a von Mises distribution with mean  $\pi$ and concentration parameter $3$ and the variables $s_i$ are uniform in $(-2,2)$, i.e. $t_i$ have support in the cylinder  with radius 1 and height between $(-2,2)$.

In all cases, the regression parameter $\bbeta$ was taken equal $5$ and the errors $\varepsilon_i$ and $\eta_i$ are i.i.d. normal with mean $0$ and standard deviation $1$.   In the smoothing procedure, the kernel was taken as the quadratic kernel $K(t)=( {15}/{16}) (1-t^2)^2 I(|x|<1)$ and we choose the bandwidth using a cross validation procedure described in Section \ref{cv}. The distance $d_g$  for these manifolds can be found  in \cite{hriv} and \cite{hrnp} and the volume density function in Section \ref{preli}. Table \ref{simu}.1 give the mean, standard deviations, mean square error for the regression estimates of $\beta$ and  the mean of the mean square error of the regression function $g$ over the 1000 replications.
\begin{center}
\begin{tabular}{ccccc}
\hline
& mean($\wbeta$) & sd($\wbeta$)& $\MSE(\wbeta)$& $\MSE(\wg)$\\
sphere case& 5.0243 & 0.0762 & 0.0064& 0.081 \\
cylinder case & 4.9845& 0.0078&  0.0003& 0.1001\\\hline
\end{tabular}
\end{center}
\vspace{-0.25cm}
\footnotesize Table \ref{simu}.1:  Performance of $\wbeta$ and $\wg$ for both models.\normalsize

In Table \ref{simu}.1 we can see a good behavior of the estimators in the two considered schemes. In all cases, the mean of the mean square error of the parametric and nonparametric estimators are small and reflect a good performance of the proposed estimators.

\subsection{Application to real data}\label{real}

In this Subsection, we applied a partially linear model to an  enviroment  dataset in order to study the atmospheric SO$_2$ pollution incidents. The variables included in the study are the direction and the speed of the wind, the temperature and the SO$_2$ concentration  in  the meteorologic station at  Villalba (Lugo in Galicia, Spain). The data was recorded daily  in each minute during the year 2009. The complete dataset has a structure of dependence in the time. Therefore to  avoid this dependence  we was considered a 2000--row  historical matrix that was constructed as in \cite{prada2}. In a previous work \cite{prada2} applied  a partial linear models to the prediction of atmospheric SO$_2$ pollution incidents in the vicinity of the coal/oil-fired power station at As Pontes (A Coru\~na in Galicia, Spain). But in this case they did not consider the direction of the wind as a directional variable.
The variables that we considered in the model was
\begin{center}
\begin{tabular}{cl}\hline
$y_i$& SO$_2$ emission is measured in $\mu g/m^3$\\
$x_{1i}$&  SO$_2$ emission in the instant $i-30$\\
$x_{2i}$&  SO$_2$ emission diference between the instant $i-35$ and  $i-30$\\
$x_{3i}$& the temperature in $^{\circ}$C\\
$t_{1i}$& wind direction in radians from the north\\
$t_{2i}$& wind speed in $m/s$\\\hline
\end{tabular}
\end{center}
\vspace{-0.25cm}
\footnotesize Table \ref{real}.1:  Enviromente variables considered in the model.\normalsize

Note that the variables $t_i=(t_{1i},t_{2i})$ have support in the cylinder. The maximum of the wind speed in this cases is $7.7$ then we consider that the variable $t$ belongs in the cylinder of high between  0 and 10. Therefore, we modeled the response variable using the following model $y_i=\beta_1 x_{1i}+\beta_2 x_{2i}+\beta_3 x_{3i}+g(t_{i})+\varepsilon_i.$

In the smoothing procedure, we considered  the quadratic kernel  and we choose the bandwidth using a cross validation procedure. Because of the computational burden of the cross-validation method, and because there is
really no need to use this method with a sample as large as 2000,  we also determined $h$ by the split sample method, i.e. by dividing the historical matrix into a 1000-member training set with odd index and a 1000-member validation set with even index, and taking for $h$
the value minimizing
\begin{eqnarray*}
SV(h) =\sum_{i=1}^{[n/2]}[(y_{2i}-\wfi_{0,E,h}(t_{2i}))-(\bX_{2i}-\wbfi_{E,h}(t_{2i}))^{\tras}\widetilde{\bbeta}]^2.
\end{eqnarray*}
where $\wbfi_{E,h}(t)=(\wfi_{1,E,h}(t),\dots,\wfi_{p,E,h}(t))$  and   $\wfi_{0,E,h}(t)$ denote the nonparametric estimators computed with bandwidth $h$ using the data with even index and $\widetilde{\bbeta}$
minimize $\sum_{i=1}^{[n/2]}[(y_{2i}-\wfi_{0,E,h}(t_{2i}))-(\bX_{2i}-\wbfi_{E,h}(t_{2i}))^{\tras}{\bbeta}]^2$ in $\bbeta$. In this case the selected bandwidth was $h_{sv}=2.5$. Table \ref{real}.2 reports the estimates values of the regression parameters and the mean and standard deviation of nonparametric estimator $\wg$ of $g$. Figure \ref{real}.1.a) shows the estimate of the regression function over a grid of 1200 points in the cylinder.
To evaluate the performance of the partial linear model, we consider a nonparametric model to explain $y_i$ based only in the variables $x_{1i}$ and $x_{2i}$ trough an unknown function $\eta$ . In this case we estimate with the Naradaya-Watson estimator with quadratic kernel. We compare the prediction error for both models computing, in the case of the full nonparametric model, we compute $EP(h) =\sum_{i=1}^{[n/2]}[(y_{2i}-\weta(x_{1,2i},x_{2,2i}))]^2$ for a grid of 100 equispaces bandwidth between 0.1 and 10. For the partial linear model we compute the $SV(h)$ for the same grid of bandwidth. As we can see in  Figure \ref{real}.1. b), the partial linear model has a better level predictive and is more stable trough the bandwidth than the full nonparametric model.

\begin{center}
\begin{tabular}{ccccc}\hline
$\wbeta_1$ & $\wbeta_2$ &$\wbeta_3$& Mean$(\widehat{g})$ & SD$(\widehat{g})$ \\
  0.9728 &  1.090 & -0.0013 & 0.1141& 0.0145\\\hline
\end{tabular}
\end{center}
\vspace{-0.25cm}
\footnotesize Table \ref{real}.2:  Estimates of regression parameter. \normalsize
\begin{center}
\hspace{-1cm}\includegraphics[scale=0.55]{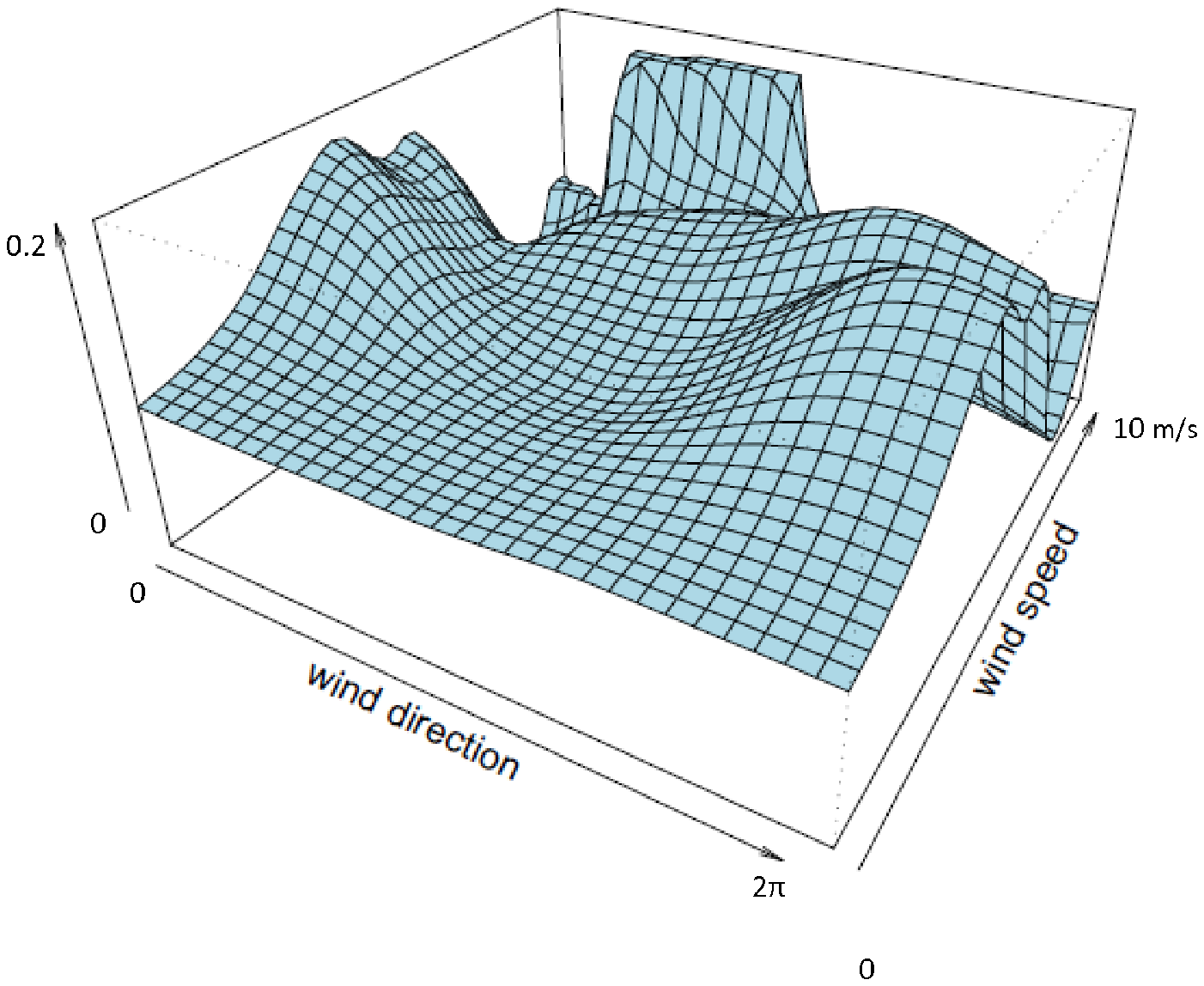}\hspace{1.4cm}\includegraphics[scale=0.45]{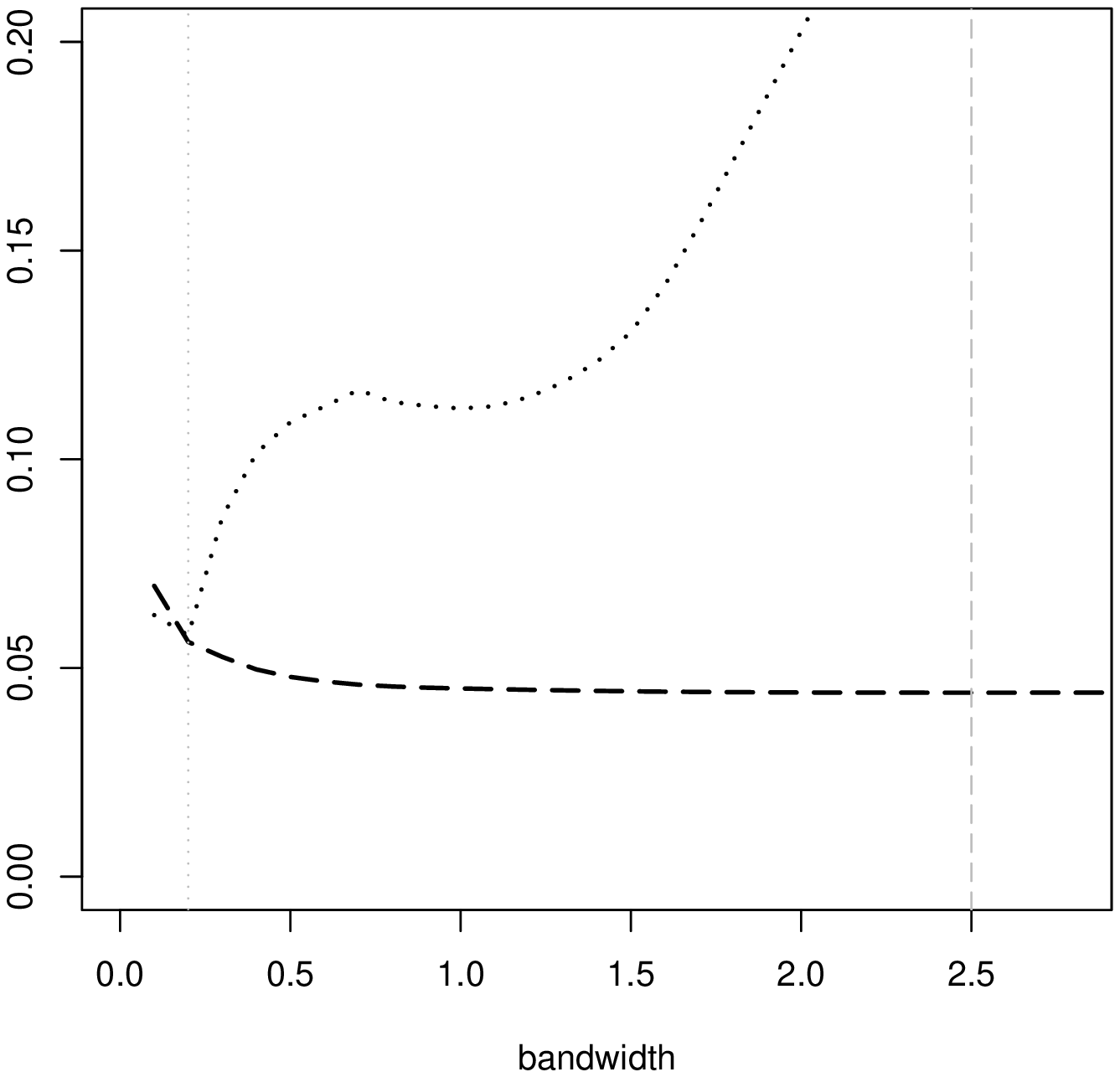}\\
\vspace{-0.8cm}
\hspace{-8cm} a) \hspace{7cm} b)
\end{center}
\vspace{-0.2cm}
\footnotesize   Figure \ref{real}.1: a) Estimates of the regression functions.  b)  Comparative of the errors: the dotted line corresponds to the full nonparametric model and the dashed line to the partially linear model. The vertical lines corresponds to the optimal bandwidth in each cases.\normalsize

\section*{Acknowledgments}

We would like to thank Mar\'\i a Leyenda Rodr\'\i guez  for the preparation of the dataset.
This work was done when the second and third authors were visiting the University of Santiago the Compostela, they are very grateful to all statistics group for their kind hospitality.  This research was partially supported by Grants X-018 from the Universidad de Buenos Aires, \sc pip \rm 1122008010216  from \textsc{conicet} and \sc pict \rm -00821  from \textsc{anpcyt}, Argentina, and also
by  Spanish Grants MTM2008-03010/MTM of the  Ministerio Espa\~nol de Ciencia e Innovaci\'on and XUGA Grants  PGIDIT07PXIB207031PR.

\thispagestyle{empty} %\cleardoublepage
\appendix
\section{Appendix}\label{proofs}

\noi\sl Lemma \ref{proofs}.1: Let $\wtfi_j(t)=\phi_j(t)-\sum_{i=1}^b w_{n,h}(t,t_i) x_{ij}$  for $1\leq j\leq p$ and  $\wtfi_0(t)=\phi_0(t)-\sum_{i=1}^b w_{n,h}(t,t_i) y_{i}$.   Under $H1$ to $H4$  we have that $\dst\max_{1\leq i\leq n}|\wtgama(t_i)|=O(h^2)+O\left(\sqrt{{\log n}/{nh^d}}\right)$ a.s.  where   $\wtgama\in\{\wtfi_j; \quad 0\leq j\leq p \}$. \rm

\vspace{0.1cm}
\noi\sl Proof of Lemma \ref{proofs}.1: \rm Let $\gamma\in\{\phi_j; \quad 0\leq j\leq p \}$ and denote by $\wgama$ the corresponding nonparametric estimator, the $\wtgama(t)=\gamma(t)-\wgama(t)$. Using analogous arguments that those considered in \cite{hriv} we have that, $\dst\sup_{t\in M_0}|E(\wtgama(t))|=O(h^2).$ Let $s_n=n^2h^{2d}\dst\sup_{t\in M_0}\left|\var(\wgama(t)\wf_n(t))\right|$
with  $\wf_n(t)={(nh^d)^{-1}}\sum_{k=1}^n{\theta^{-1}_t(t_k)}K(d_g(t,t_k)/h)$, by results obtained in \cite{hrnp} we have that $s_n=O(nh^d)$.
By $H1$, $\inf_{t\in M}\frac{1}{h^d}E\left(\frac{1}{\theta_t(t_1)}K(d_g(t,t_1)/h)\right)\geq A> 0.$ Then, it follows in analogous way  that the proof of  Lemma 3.1 in  \cite{fv}.\square

\vspace{0.25cm}

\noi\sl Lemma \ref{proofs}.2: Under $H1$ to $H4$  we have that $n^{-1}\wtX^{\tras}\wtX\convprob \Sigma$.   \rm

\vspace{0.1cm}
\noi\sl Proof of Lemma \ref{proofs}.2: \rm The   element $l,s$  of $n^{-1}\wtX'\wtX$ can be written as
\begin{eqnarray*}
(n^{-1}\wtX^{\tras}\wtX)_{ls}=n^{-1}\wtX_l^{\tras}\wtX_s=n^{-1}\left(\sum_{i=1}^n\eta_{il}\eta_{is}+\sum_{i=1}^n\wtfi_l(t_i)\eta_{is}+\sum_{i=1}^n\wtfi_s(t_i)\eta_{il}+\sum_{i=1}^n\wtfi_l(t_i)\wtfi_s(t_i)\right)\\
\end{eqnarray*}
where $\wtfi_j(t)=\phi_j(t)-\wfi_j(t)$.  We need to show that all terms except the first term converge to zero and by applying the strong law of large numbers we get that $n^{-1}\sum_{i=1}^n\eta_{il}\eta_{is}\convprob \Sigma_{ls}$. Since Lemma \ref{proofs}.1 and  the fact that $n^{-1}\sum_{i=1}^n\eta^2_{il}\convprob \Sigma_{ll}$ and  using the Cauchy-Schwarz inequality we get the result.\square

\vspace{0.25cm}

\noi\sl Lemma \ref{proofs}.3: Under  $H1$ to $H3$,  we have that $\max_{1\leq i,j\leq n }|w_{n,h}(t_i,t_j)|=O((nh^d)^{-1})$.   \rm

\vspace{0.1cm}
\noi\sl Proof of Lemma \ref{proofs}.3: \rm
Using the results obtained in \cite{hriv} and \cite{hrnp} we have that
\begin{eqnarray}
\sup_{t\in M}\left|\frac{1}{nh^d}\sum_{i=1}^n \dst\frac{1}{\theta_t(t_i)}K(d_g(t,t_i)/h) -\frac{1}{h^d}E\left(\frac{1}{\theta_t(t_1)}K(d_g(t,t_1)/h)\right)\right|=o(1)\quad \mbox{ a.s.}\label{pesos1}\\
\inf_{t\in M}\frac{1}{h^d}E\left(\frac{1}{\theta_t(t_1)}K(d_g(t,t_1)/h)\right)\geq A> 0.\label{pesos2}
\end{eqnarray}
Then by (\ref{pesos1}) and (\ref{pesos2}) and the boundedness of $K$ and $\theta_t$, the lemma holds.\square

\vspace{0.25cm}
\noi\sl Remark \ref{proofs}.4: \rm Note that by Lemmas \ref{proofs}.1 and \ref{proofs}.3 and using Lemma A.1 in \cite{liang}; we have that
 $\dst\max_{1\leq i\leq n}|\gamma(t_i)-\sum_{k=1}^nw_{n,h}(t_i,t_k)\gamma(t_k)|=O(h^2)+O\left(\sqrt{{\log n}/{nh^d}}\right) $ a.s. for any $\gamma\in\{\phi_j; \quad 0\leq j\leq p \}$.

\vspace{0.25cm}

\noi\sl Proof \ref{asymp}.1: \rm
We can write $
\sqrt{n}(\wbeta-\bbeta)=(n^{-1}\wtX^{\tras}\wtX)^{-1}n^{-1/2}\left[ A_{1n}-A_{2n}+A_{3n}\right]$ where
$$
A_{1n}=\sum_{i=1}^n\wtX_ig^*(t_i)\quad A_{2n}=\sum_{i=1}^n\wtX_i\left(\sum_{i=1}^nw_{n,h}(t_i,t_j)\varepsilon_j\right) \quad A_{3n}=\sum_{i=1}^n\wtX_i\varepsilon_i
$$
and $g^*(t)=g(t)-\sum_{i=1}^nw_{n,h}(t,t_i)g(t_i)$. Using Lemmas \ref{proofs}.1 to \ref{proofs}.3, the asymptotic behavior of $A_{1n}, A_{2n}$ and $A_{3n}$ can be obtained in the same way that in \cite{av}. Specifically, considering the assumptions imposed on $h$,  we can obtained that
\begin{eqnarray*}
A_{n1}&=&O(nh^4+h^{-d}\log^2 n)+O(n^{1/2}h^2\log n+h^{-d/2}\log^2 n)+O(n^{1/2}h^2h^{-d/2}\log n)+O(h^{-d}\log^2 n)=o(n^{1/2})\\\
A_{n2}&=&O(n^{1/2}h^2h^{-d/2}\log n+h^{-d}\log^2 n)+O(h^{-d/2}\log^2 n)+O(h^{-d}\log^2 n)=o(n^{1/2})\\
A_{n3}&=&O(n^{1/2}h^2\log n)+O(h^{-d/2}\log^2 n)+\sum_{i=1}^n\eta_i\varepsilon_i+O(h^{-d/2}\log^2 n)=\sum_{i=1}^n\eta_i\varepsilon_i+o(n^{1/2})\\
\end{eqnarray*}
Finally, the central limit theorem gives the desired result. \square

\end{document}